\input amstex.tex

\input amsppt.sty

\TagsAsMath

\magnification=1200

\hsize=5.0in\vsize=7.0in

\hoffset=0.2in\voffset=0cm

\nonstopmode

\document

\def\R{ \Bbb R}

\input amstex.tex
\input amsppt.sty
\TagsAsMath \NoRunningHeads \magnification=1200
\hsize=5.0in\vsize=7.0in \hoffset=0.2in\voffset=0cm \nonstopmode

\document

\topmatter

\title{An invariant set  in energy space  for supercritical NLS in 1D}
\endtitle

\author
Scipio Cuccagna
\endauthor

\address
DISMI University of Modena and Reggio Emilia, via Amendola 2,
Padiglione Morselli, Reggio Emilia 42100 Italy\endaddress \email
cuccagna.scipio\@unimore.it \endemail

\abstract We consider radial  solutions of a mass supercritical
monic NLS and we    prove the existence of a set, which looks like a
hypersurface, in the space of   finite energy functions, invariant
for the flow   and formed by solutions which converge to ground
states.

\endabstract

\endtopmatter

 \head \S 1 Introduction \endhead

We consider a monic supercritical NLS
$$iu_t +  u_{xx} +  |u| ^{p-1} u=0\, ,\, (t,x)
\in \Bbb R \times \Bbb R \, , \,  u(t,x)\equiv u(t,-x) \, ,
\,5<p<\infty .\tag 1.1$$ We ignore translation and   consider only
even solutions $u(t,x)\equiv u(t,-x)$ of (1.1):
  by $H^1_r (\Bbb R,\Bbb C)$
  we will mean  the space of finite  energy even
functions. (1.1) admits ground states solutions $    e^{
   i t\omega
+i\gamma  }  \phi _{\omega }  (x ) $, with $\phi _\omega (x)= \omega
^{\frac 1{2(p-1)}}  {(\frac {p+1}2 )^{\frac 1{p-1}}}{ \text{sech}\,
^{\frac 2{p-1}}(\frac{p-1}2\sqrt{\omega} x)}.$   Let
$$\aligned & G= \big \{ e^{i \gamma  }  \phi _{\omega }  (x ) : \, \omega >0; \,  \gamma
 \in \Bbb R  \big \} \subset
H^1_{r}(\Bbb R,\Bbb C).
\endaligned \tag 1.2$$
  For any initial datum $u(0,x)\in H^1_r (\Bbb R,\Bbb C)$
close   to $G$, for some time the corresponding solution $u(t,x)$
remains close to $G$ and can be written in a canonical way as a
varying  ground state plus a reminder term:
$$\aligned & u(t,x) = e^{i \int _0^t\omega (s) ds
+i\gamma (t)} (\phi _{\omega (t)} (x)+ r(t,x))   .
\endaligned \tag 1.3$$
The orbits in $G$ are  unstable and $u(t,x)$ can
  blow up in finite time  \cite{BC},  so   (1.3) in general does not
persist for all $t$. We will prove:

\proclaim{Theorem 1.1} There exist   a $X\subset H^1 _r(\Bbb R ,\Bbb
C)$   such that:

{\item{$\bullet $} } $G\subset X$;

{\item{$\bullet $} }  $X$ is invariant by the flow;

{\item{$\bullet $} } $X$ looks like a hypersurface, in the following
sense: for any $ g_0\in G$   there exists a neighborhood $U$ of $
g_0 $ in $H^1 _r(\Bbb R,\Bbb C)$ such that there is
$\widetilde{X}\subseteq X\cap U$ with $ \widetilde{X}$   the graph
of a real valued   function, non necessarily continuous, defined on
a real closed hyperplane through $g_0$ in $H^1_r (\Bbb R,\Bbb C)$;

{\item{$\bullet $} } For any $ g_0=e^{i \gamma _{0} }  \phi _{\omega
_{0}}  (x ) \in G$ there are $ C>0$ and $\epsilon _0>0$, which
depend  only on $\omega _0$, such that for any $0<\epsilon <\epsilon
_0 $ if we pick $u_0\in X$ with $\| u_0-g_0\| _{H^1  (\Bbb
R)}<\epsilon $ then the corresponding solution $u(t,x)$ is globally
defined and contained in $X $, can be written in a canonical way in
the form (1.3), and  we have
$$\align & \| r(t  )\| _{H^1(\Bbb R,\Bbb C)}+ \| r \| _{L^4_t(\Bbb R_+,L^\infty  _x(\Bbb R,\Bbb C))}+ |(\omega  _0,\gamma _0) -(\omega  (t),\gamma (t))
|<C\epsilon  .\tag 1\endalign
$$   The   limit
$$\lim _{t\to +\infty }(\omega  (t),\gamma (t))=(\omega  _\infty ,\gamma _\infty )   \tag 2 $$
exists and there exists $r_\infty \in H^1_r (\Bbb R,\Bbb C)$ with $
\| r_ \infty \| _{H^1(\Bbb R,\Bbb C)}<C\epsilon  $    such that
$$\lim _{t\to +\infty }  \| e^{ i \int
_0^t \omega   (\tau ) d\tau  +i\gamma (t)}r(t)  -e^{ it
\partial ^2_x } r_\infty\| _{H^1 (\Bbb R, \Bbb C  )}  =0. \tag 3$$

\endproclaim
{\it Remark.} In the subspace  of $H^1_r\times H^1_r$ formed by
pairs $(u , \overline{u})$, the hyperplane at $ g_0=e^{i \gamma _{0}
} \phi _{\omega _{0}}  (x ) $ is spanned  by  $N_g(H_{\omega
_0})\oplus \Bbb R \sigma _1\xi (\omega _0) \oplus L^2_c(H_{\omega
_0})$, with the various terms introduced in  \S 2.

\medskip {\it Remark.} We emphasize that all the functions considered
in this paper are even in $x$.

\bigskip

\noindent Theorem 1.1 is related   Tsai \& Yau \cite{TY}, Schlag
\cite{S} and Krieger \& Schlag \cite{KS}.   \cite{KS}
  for (1.1) proves the existence of a Lipchiz
hypersurface of initial data  $u_0$ with    $\langle x \rangle u_0
\in H^1(\Bbb R ) \cap W^{1,1}(\Bbb R ) <\infty $, such that the
corresponding solutions $u(t,x)$ converge  to ground states. The
stronger decay hypothesis on the initial data allows to control the
rate of convergence of $\omega (t)$ to its limit, and also the rate
of convergence of the motion of the ground state to the inertial
asymptotic motion. For data  $H^1(\Bbb R )$ or in the smaller space
$H^1_r(\Bbb R )$ the method in \cite{KS} does not work. We consider
only even initial data to eliminate spatial motion of the ground
state. So the velocity is zero and we trivialize one of the
difficulties. The problem with $\omega (t)$ however remains. We
obtain our result by means of Schauder fixed point theorem applied
to an appropriate functional. Unfortunately, due to the fact that
$u_0\in H^1_r(\Bbb R)$ and to the lack of sufficient  control on
$\omega (t)$, we are not able to show that the  functional is a
contraction, which would yield $X=\widetilde{X}$  and  some
regularity for the hypersurface. It would be nice to prove that $X$
is a continuous hypersurface, and then, given a small ball $B\subset
H^1_r(\Bbb R )$ of center $g\in G$, to study the behavior of
solutions which start in $B\backslash X$. During the review process
of this paper we learned of the work by  Beceanu \cite{B} which
proves an analogous result to the present one for solutions $u(t)\in
H^1(\Bbb R^3)\cap L^{2,1}(\Bbb R^3)$, in the notation below, for the
cubic NLS treated in \cite{S}. The result in \cite{B} is stronger
than ours in two respects: $X$ is indeed a Lipschitz hypersurface in
$H^1(\Bbb R^3)\cap L^{2,1}(\Bbb R^3)$, and there is no requirement
of spherical or other symmetries. The proof in \cite{B} does not
work in our 1 dimensional setting  for solutions $u(t)\in H^1(\Bbb R
)\cap L^{2,1}(\Bbb R )$. We remark that the endpoint Strichartz
estimate needed in \cite{B} is a corollary of the transposition to
linearizations of the NLS of the following material:
  Yajima's $L^p$ theory of wave operators \cite{Y1,Y2} transposed in \cite{C4,C5};
Kato smoothness theory \cite{K}, applied in   Proposition 4.1
\cite{CPV}. Furthermore, in cases when they cannot be derived
directly from bounds on wave operators, as for example Lemma 3.1
below, Strichartz estimates for the linearization $H_\omega$ in
(2.2) can be proved with a standard $TT^*$ argument, using an
appropriate bilinear form, see the proof of Lemma 3.1 in
\cite{C1,C3} . For other results related to the present paper see
\cite{Co,Ma} and references therein.

 In the last section we
list a series of errata in paper \cite{C1}. In particular the
present paper is based on \cite{C3}, which is  a thorough revision
of \cite{C1}.

 We write
$R_H(z)=(H-z)^{-1}$ and
   $\langle x \rangle =(1+|x|^2) ^{\frac
12}$.
    We set $\|
u\|_{H^{k,s}}:= \|\langle x \rangle  u\|_{H^{k }} $. We set $L^{2,s}
=H^{0,s}.$
 We set $
\langle f ,
  g\rangle =\int _{\Bbb R}{^tf(x)} {g(x)} dx,$  with $f(x)$
  and $g(x)$ column vectors and with $^tA$ the transpose.
$  W ^{1,p}(\Bbb R)$ is the set of tempered  distributions $f(x)$
with derivative $f(x),f'(x)\in L ^{ p}(\Bbb R)$. $\Cal W^{k,p }(\Bbb
R)$ is the space of tempered distributions $f(x)$ such that $
(1-\partial _x^2) ^{k/2}f\in L^p (\Bbb R)$. Recall that $  W
^{1,p}(\Bbb R)=\Cal W^{k,p }(\Bbb R)$ exactly for $1<p<\infty .$
$\dot W ^{1,p}(\Bbb R)$ is the set of tempered  distributions $f(x)$
with derivative $f'(x)\in L ^{ p}(\Bbb R)$.

\head \S 2 Linearization and spectral decomposition \endhead

We plug the ansatz  (1.3) in (1.1) obtaining, for $n(  r,\overline{
r })=O(r^2),$

$$\aligned &
  i r_t  =-  r _{xx} +\omega (t) r -
    \frac{p+1}{2}\phi _{\omega (t)}  ^{p-1}
  r  (t,y) -\frac{p-1}{2}\phi _{\omega (t)}  ^{p-1} \overline{ r
  }
  + \\& +
\dot \gamma(t)  \left ( \phi _{\omega (t)} +r  \right )   -i\dot
\omega (t)
\partial _\omega \phi   _{\omega (t)}
 +n(  r,\overline{ r })
\endaligned   \tag 2.1 $$
Let $\sigma _1=\left[ \matrix 0 &
1  \\
1 & 0
 \endmatrix \right] \, ,
\sigma _2=\left[ \matrix  0 &
i  \\
-i & 0
 \endmatrix \right] \, ,
\sigma _3=\left[ \matrix  1 &
0  \\
0 & -1
 \endmatrix \right]   .$
The linearization is
$$ \aligned & H_\omega  =\sigma _3
\left ( - {d^2}/{ dx^2}+\omega \right )+\omega V(\sqrt{\omega} x)\\&
V(x)=-  \left (\sigma _3
   (p+1)  -i \sigma _2 ({p-1})\right )
  (p+1)    2^{-2}{ \text{sech}\,
^{2}(\frac{p-1}2  x)} . \endaligned\tag 2.2$$ By (2.1), for
 $^tR
 =(  r,
\bar r )$, $^t  \Phi  = (\phi _{\omega  }, \phi _{\omega } )$ and
$^t   N (R )  = (n(  r,\overline{ r }), -\overline{n( r,\overline{ r
})} )$,

$$\aligned & i  R _t  =H_{\omega  (t)}   R
+\sigma _3  \dot \gamma R      +\sigma _3   \dot \gamma \Phi  - i
\dot \omega \partial _\omega \Phi +  N (R ).\endaligned \tag 2.3
$$

By implicit function theorem we impose  $R(t )\in N_g^\perp
(H_{\omega (t) }^\ast )$, with $N_g$ the generalized kernel.   We
state the following known result:

\proclaim{Theorem 2.1} Let us consider the operator  $H_\omega $ in
(2.2) acting on $R\in L^2(\Bbb R,\Bbb C^2)$:

{\item {(1)}} The continuous spectrum of $H_\omega $ is $\Bbb
R\backslash (-\omega ,\omega )$.     0 is an eigenvalue and there
are   two simple eigenvalues $\pm i \mu (\omega )$, with $ \mu
(\omega )>0.$

  {\item {(2)}}  $N_g(H_{\omega } )$ is spanned by
 $ \{ \sigma _3
\Phi _\omega   ,  \partial _\omega \Phi  _\omega     , \partial _x
\Phi _{\omega  }   , \sigma _3x\Phi _\omega      \}   .$

{\item {(3)}} $\pm \omega    $ are not resonances and  $\{ 0, i \mu
(\omega ),-i \mu (\omega )\}$  are the only eigenvalues.

\endproclaim
 For (1) and (2) see \cite{W}, for (3) see \cite{KS}.
  Let  $\xi  (\omega  ,x)$ be an eigenvector of $  i\mu  (\omega )$.
  Notice that $\mu  (\omega )=\omega \mu  (1 )$.
  Recalling that $ \langle f ,
  g\rangle =\int _{\Bbb R}{^tf(x)} {g(x)} dx,$ we have:

\proclaim{Lemma 2.2}
 The
 eigenvector  $\xi  (\omega  ,x)$ can be chosen so that
$\langle \xi (\omega ),\sigma _3\xi (\omega )\rangle = i\lambda _1 $
with $\lambda _1\in \Bbb R  \backslash \{ 0 \}  $ a fixed number.
The function $(\omega , x )  \to \xi (\omega , x)$ is $C^2$ and
$|\xi (\omega , x)| < c\sqrt{\omega} e^{-a\sqrt{\omega }|x|}$ for
fixed $c>0$ and $a>0$.
     $\sigma _1\xi (\omega ,x) =\overline{\xi (\omega ,x)}$
      generates $\ker (H
 _\omega
+i\mu (\omega ))$ with $\langle \sigma _1\xi ,\sigma _3\sigma _1\xi
\rangle = -i \lambda _1$.  We have
 $H_{\omega }$   invariant
  decompositions
$$ L^2(\Bbb R,\Bbb C^2)= L^2_d( \omega ) \oplus  L_c^2( {\omega })
\text{ and } L^2(\Bbb R,\Bbb C^2)=N_g(H_{\omega })\oplus
 N_g^\perp (H_{\omega  }^\ast  )\tag 1
$$ with $ L^2_d( \omega ) =N_g(H_{\omega }) \oplus \big ( \oplus _{
\pm }\ker (H _\omega \mp i\mu (\omega  )) \big )$ and $L_c^2(
{\omega })=\left [\sigma _3L_d^2(\omega )\right ]^\perp   $.
\endproclaim
{\it Proof.} The   decomposition (1) is a consequence of Theorem
2.1. Let $\xi (x) $ be a generator of  $\ker (H
 _1
-i\mu (1 ))$.  Since  both $\overline{\xi} (x)$ and $\sigma _1\xi
(x)\in \ker (H
 _1
+i\mu (1 ))$, we can normalize $\xi (x)$ so that $\overline{\xi}
(x)=\sigma _1\xi (x)$. Then $^t\xi (x)=(v(x),\overline{v}(x)).$ Then
$\langle \xi ,\sigma _3\xi \rangle =\int (v^2-\overline{v}^2)dx =i
\lambda _1  $ with $\lambda _1\in \Bbb R  \backslash \{ 0 \}  $.
Notice that $\lambda _1\neq 0$ since otherwise $\langle \xi ,\sigma
_3f \rangle =0$ for any $f$ would follow from the fact that $\langle
\xi ,\sigma _3f \rangle =0$ for any $f\in N_g(H _1)\oplus  L_c^2(
1)$ and for $f=\sigma _1\xi .$
  Finally set   $\xi
(\omega  ,x)=\sqrt{\omega} \xi  (1 ,\sqrt{\omega}x).$ The rest is
standard.

\bigskip
We   denote by $P_d(\omega )$ (resp. $P_c(\omega )$) the projection
on $L^2_d( \omega )$ (resp. $L^2_c( \omega )) $ associated to the
splitting in (1) Lemma 2.2. By $ N_g (H_{\omega }^\ast )=\sigma
_3N_g (H_{\omega} ) $, the condition $R(t,x)\in N_g^\perp (H_{\omega
(t)}^\ast )$  and (2.3) imply    the modulation equations:

$$\aligned &
i\dot \omega \,  d(\|  \phi _\omega \| _2^2)/d\omega  =i\dot \omega
\big \langle R ,\partial _\omega \Phi _\omega  \big \rangle +\big
\langle \sigma _3 \dot \gamma R   +  N (R ) , \Phi _\omega
 \big \rangle
\\&  \dot \gamma \,  d(\|  \phi _\omega \| _2^2)/d\omega  =i\dot \omega \big
\langle R , \sigma _3\partial _\omega  ^2\Phi _\omega  \big \rangle
- \big \langle \sigma _3  \dot \gamma R +   N (R ) , \sigma
_3\partial _\omega  \Phi _\omega  \big \rangle
  .
\endaligned    $$
 By elementary computations, see \cite{C2}, there are real valued exponentially decreasing
functions $\alpha (\omega ,x)$ and $\beta (\omega ,x)$ such that

$$\align &
\Cal M (\omega ,R)     \left [ \matrix  i\dot \omega \\ -\dot \gamma
\endmatrix  \right ]  = \left [ \matrix   \langle n(r,\overline{r}) -n( \overline{r},r), \phi _\omega \rangle  \\
  \langle n(r,\overline{r})+n( \overline{r},r),   \partial _\omega \phi _\omega
\rangle
\endmatrix \right ] \text{ with}\\& \Cal M (\omega ,R)=
  d(\| \phi _\omega \| ^2_2)/d\omega   +\left [    \matrix
\langle r+\overline{r}, \alpha (\omega ) \rangle
&   \langle r-\overline{r} , \phi _\omega \rangle  \\
 \langle r-\overline{r}, \beta (\omega ) \rangle  & \langle
r+\overline{r} , \partial _\omega \phi _\omega \rangle
\endmatrix \right ]   .\tag 2.4
\endalign   $$
Since in the sequel we deal  with $R(t)$   such that $\| R \|
_{L^\infty _tL^2_x} $ is  small and such that $\omega $ remains in
a bounded domain, we get
$$\aligned & \left [ \matrix  i\dot \omega (t) \\ -\dot \gamma (t)
\endmatrix  \right ]  = \left [ \matrix
 i\dot {\widetilde{\omega}} (\Cal R) \\
 -\dot {\widetilde{\gamma}} (\Cal R)
\endmatrix  \right ] \text{ with}\\&
\left [ \matrix
 i\dot {\widetilde{\omega}} (\Cal R) \\
 -\dot {\widetilde{\gamma}} (\Cal R)
\endmatrix  \right ] :=
M (\omega ,R)\left [ \matrix   \langle n(r,\overline{r}) -n( \overline{r},r), \phi _\omega \rangle  \\
  \langle n(r,\overline{r})+n( \overline{r},r),   \partial _\omega \phi _\omega
\rangle
\endmatrix \right ]
\text{ with }  M (\omega ,R):=\\&   =\Cal M ^{-1}(\omega ,R) =
 \left (  d(\| \phi _\omega \| ^2_2)/d\omega  \right ) ^{-1}
 \left ( 1+O (\| R \|
_{L^\infty _tL^2_x}) +O(\| \omega  -\omega _0\| _{L^\infty
_t})\right ) .
\endaligned   \tag 2.5$$

\proclaim{Lemma 2.3} We can write  $   R(t)=  f (t) +\zeta (t) $
with $f(t)\in L^2_c( { \omega (t)})$ and $\zeta (t,x) =z_{+}(t)  \xi
(\omega (t),x) +   z _{-}(t) \sigma _1 \xi (\omega (t), x)$. $ R =
\sigma _1\overline{ R } $ implies  $z_\pm (t)\in \Bbb R$   and $ f =
\sigma _1\overline{ f  } $.

\endproclaim
{\it Proof. }By $R(t,x)\in N_g^\perp (H_{\omega (t)}^\ast )$ and
setting $f=P_c(\omega (t))R$ we get $   R(t)=  f (t) +\zeta (t) $
for an $\zeta (t ) =z_{+}(t)  \xi (\omega (t) ) +   z _{-}(t) \sigma
_1 \xi (\omega (t) )$.   $z_\pm (t)\in \Bbb R$ and  $ f  = \sigma
_1\overline{ f  } $ follow  by $ \xi = \sigma _1\overline{ \xi } $,
$ \sigma _1\xi = \overline{ \xi } $, $\sigma _1 L^2_c(\omega (t))=
L^2_c(\omega (t)), $   $ \overline{L^2_c(\omega (t))}= L^2_c(\omega
(t))  $  and
$$z_{+}  \xi   +   z _{-}  \sigma _1 \xi
+f=R=\sigma _1\overline{R}=\overline{ z_{+}} \sigma _1\overline{ \xi
} +\overline{ z _{-} }  \overline{\xi} +\sigma _1\overline{ f}.$$

\bigskip We have from (2.3) and Lemma 2.3
$$\aligned & if_t=  H_{\omega(t)}f+\sigma_3 \dot \gamma R +N(R)+   \sigma _3\dot \gamma \Phi_{\omega(t)}- i\dot\omega
\partial_\omega \Phi_{\omega (t)}
\\ & + i(z_+\mu (\omega(t))- \dot z_+)\xi(\omega(t))
-i( {z}_-\mu (\omega(t))+ \dot{ z}_-)\sigma_1\xi(\omega(t))
\\ &   -i \dot \omega (z \partial _\omega
\xi(\omega(t)) + \bar z \sigma _1 \partial _\omega \xi(\omega(t)) ).
\endaligned$$
We apply $\langle \cdot , \sigma _3\xi \rangle $  and $\langle \cdot
, \sigma _3\sigma _1\xi \rangle $. Setting $ d_1= -\lambda _1^{-1}$
with $\lambda _1$ the constant in Lemma 2.2, we get the discrete
mode equations:

$$\aligned &  \dot z _{\pm }(t)\mp  \mu  (\omega (t)) z _{\pm }(t) =
  d_1 i \dot \omega \langle f(t ), \sigma _3 \sigma _1^{\frac{1\mp
1}{2} }\partial _\omega \xi (\omega (t) ) \rangle +
\\& d_1\langle \sigma_3\dot \gamma R+  N(R)- i\dot \omega (t)\left [
z_{+}(t)
  +   z _{-}(t) \sigma _1
\right ] \partial _\omega  \xi (\omega (t) ), \sigma _3\sigma
_1^{\frac{1\mp 1}{2}  } \xi (\omega (t) )   \rangle     .
  \endaligned \tag 2.6$$
We fix an $\omega _0$. Setting $\omega =\omega (t)$ and $\ell
(t)=\omega (t) -\omega _0 +\dot
  \gamma (t)$
we get

$$\aligned & \left [ i\partial _t  - \left ( H_{\omega _0} +\ell (t)
 P_c(\omega _0)\sigma _3
  \right ) \right ] f = P_c(\omega )\sigma _3\dot \gamma (z _{+}+
z_{-} \sigma _1  ) \xi +\Cal N(R) +\\&   +\left ( \omega _{0}
V(\sqrt{\omega _0} x)- \omega V(\sqrt{\omega} x)\right ) f+i\dot
\omega
\partial _\omega P_c(\omega )f + \ell(t)  \left ( P_c(\omega )-P_c(\omega _0) \right
) \sigma _3
 f.\endaligned \tag 2.7
$$
To correct the fact that $[P_c(\omega _0)\sigma _3, H_{\omega
_0}]\neq 0$, we   split $f\in L_c^2( {\omega (t)})$ into
$$ f=f_d+f_c\in L^2_d( {\omega _0})\oplus L^2_c( {\omega _0}).
\tag 2.8 $$ Then splitting $P_c(\omega _0)=P_+(\omega _0)+P_-(\omega
_0)$, with the two terms the projections on the positive and the
negative part of the continuous spectrum,  see Lemma 5.12 \cite{C1}
or Appendix B \cite{C3} or also \cite{BP}, we get

$$\aligned & \left [ i\partial _t  - \left ( H_{\omega _0}
+\ell (t) (P_+(\omega _0)-P_-(\omega _0))
  \right ) \right ] f _c=  P_c(\omega _0)\{ P_c(\omega )\sigma _3\dot \gamma (z _{+}+
z_{-} \sigma _1  ) \xi \\& +  N(R) + \left ( \omega _{0}
V(\sqrt{\omega _0} x)- \omega V(\sqrt{\omega} x)\right ) f+i\dot
\omega
\partial _\omega P_c(\omega )f \\& + \ell(t)  \left ( P_c(\omega )-P_c(\omega _0) \right
) \sigma _3
 f  + \ell (t)(P_+(\omega
_0)-P_-(\omega _0)-P_c(\omega _0)\sigma _3) f_c \} . \endaligned
\tag 2.9
$$
Now $[ P_+(\omega _0)-P_-(\omega _0) , H_{\omega _0}]=0$. We will
use the following elementary lemma.

\proclaim{Lemma 2.4} Fix $\alpha \in (0,1)$. Then there exists a
small $\delta (\alpha )>0$ such that for any fixed $\omega _0\in
(\alpha , 1/\alpha)$ and for any $\omega  $   with $ |\omega -\omega
_0 | \le \delta (\alpha ) $ there exist constants $C_N(\alpha )$
such the following holds: for any $f_c \in L^2_c( {\omega _0})$
there exists exactly one $f_d \in L^2_d( {\omega _0})$ so that $f
=f_c +f_d  \in L^2_c( {\omega })$ and for any $q\in [1,\infty ]$ we
have
$$ \| f_d \| _{L^q _x}\le C_N(\alpha )|\omega
-\omega _0| \, \| \langle x \rangle ^{-N}f_c \| _{L^2 _x} .\tag
2.10$$ Furthermore, if $\overline{f_c}=\sigma _1 f_c$, then we have
$\overline{f_d}=\sigma _1 f_d.$
\endproclaim
{\it Proof. } For $f_c =P_d(\omega  )f_c + P_c(\omega  )f_c $ we
seek $f_d \in L^2_d( {\omega _0})$ with
  $P_d(\omega   ) f_d
=-P_d(\omega   ) f_c $. We have $P_d(\omega   )P_d(\omega _0
)=P_d(\omega _0  )+   (P_d(\omega   )- P_d(\omega _0 )) P_d(\omega
_0 )$. Since $P_d(\omega   )- P_d(\omega _0 )=O(\omega    -
 \omega _0  )$ in any norm, we see for the ranks, $\text{Rk}\left( P_d(\omega   )P_d(\omega _0
)
 \right )=\text{Rk}\left( P_d(\omega _0
)
 \right )=\text{Rk}\left( P_d(\omega   )
 \right )$, so $P_d(\omega   )P_d(\omega _0
):L^2_d( {\omega _0})\to L^2_d( {\omega  })$ is an isomorphism and
  $f_d $ exists
 unique. Next,
 $$ -P_d(\omega   ) f_c = (P_d(\omega _0 )-P_d(\omega   ))
 f_c =f_d  +   (P_d(\omega   )- P_d(\omega _0 ))f_d $$
 implies
 $  \| f_d \| _{q} (1-C|\omega  -\omega _0|)\le \| (P_d(\omega _0 )-P_d(\omega   ))
 f_c \| _{q} \lesssim |\omega
-\omega _0| \, \| \langle x \rangle ^{-N}f_c \| _{2}.$
  Let $J$ be either
$\sigma _1 $ or the conjugation operator $Jh=\overline{h}$. Then, in
either case $ [P_d(\omega  ),J]=[P_c(\omega  ),J]=0$ for any
$\omega$. This implies $\overline{f_d}=\sigma _1 f_d $.

\head \S 3 Spacetime estimates for $H_\omega $ \endhead

We will need the following estimates, proved in \cite{C3}.

 \proclaim{Lemma 3.1 (Strichartz estimate)} Let  $\Cal W^{k,p
}(\Bbb R)$ be the space of tempered distributions $f(x)$ such that $
(1-\partial _x^2) ^{k/2}f\in L^p (\Bbb R)$. Then there exists a
positive number $C=C(\omega )$ upper semicontinuous in $\omega $
such that for any $k\in [0,2]$:
  {\item {(a)}}
 for any $f\in
L^2_c( {\omega })$, $ \|e^{-itH_{\omega }} f\|_{L_t^4\Cal
W_x^{k,\infty }\cap L_t^\infty H_x^k}\le C\|f\|_{H^k}.$
 {\item {(b)}}
  for any $g(t,x)\in
S(\R^2)$,
$$
\|\int_{0}^te^{-i(t-s)H_{\omega }} P_c( {\omega
})g(s,\cdot)ds\|_{L_t^4\Cal W_x^{k,\infty }\cap L_t^\infty H_x^k}
\le C\|g\|_{L_t^{4/3}\Cal W_x^{k,1 }+L_t^1H_x^k}.
$$\endproclaim

 \proclaim{Lemma 3.2} For any $k$ and  $\tau >3/2$  $\exists$
  $C=C(\tau ,k,\omega )$ upper
semicontinuous in $\omega $ such that:
   {\item {(a)}}
  for any $f\in S(\R)$,
$$\align &
\| e^{-itH_{\omega }}P_c(H_\omega )f\| _{L_{  t}^2 H_x^{k, -\tau}}
\le
 C\|f\|_{H ^{k}} .
\endalign $$
  {\item {(b)}}
  for any $g(t,x)\in
 {S}(\R^2)$
$$ \left\|\int_\Bbb R e^{itH_{\omega }}
P_c(H_\omega )g(t,\cdot)dt\right\|_{H^k_x} \le C\| g\|_{L_{  t}^2
H_x^{k, \tau}}.
$$\endproclaim

\proclaim{Lemma 3.3} For any $k$ and  $\tau >3/2$  $\exists$
  $C=C(\tau ,k,\omega )$ as
above such that $\forall$ $g(t,x)\in {S}(\R^2)$
$$\align &  \left\|  \int_0^t e^{-i(t-s)H_{\omega
}}P_c(H_\omega )g(s,\cdot)ds\right\|_{L_{  t}^2 H_x^{k, -\tau}} \le
C\|  g\|_{L_{  t}^2 H_x^{k, \tau}}.\endalign
$$
\endproclaim

\proclaim{Lemma 3.4} $k$ and  $\tau >3/2$  $\exists$
  $C=C(\tau ,k,\omega )$ as
above such that $\forall$ $g(t,x)\in {S}(\R^2)$
$$\align  &
\left\|\int_0^t e^{-i(t-s)H_{\omega }}P_c(H_\omega )g(s,\cdot)ds
  \right\|_{
L_t^\infty L_x^2\cap  L^{4 }_{t}(\Bbb R ,W ^{k,\infty }_x) } \le
C\|g\|_{L_t^2H_x^{k,\tau}}.
\endalign$$
\endproclaim

 \proclaim{Lemma 3.5} In Lemmas 3.1 (b),   3.3 and 3.4
the estimates continue to hold if we replace in the integral $[0,t]$
with $[t,+\infty ).$
\endproclaim

\head \S 4 Functional setting and integral formulation \endhead

From now on in the paper all the functions we consider are even in
$x$. We want to build a set $X$  of  special solutions of (1.1)
which for all times are approximate ground states $u(t,x) = e^{i
\int _0^t\omega (s) ds +i\gamma (t)} (\phi _{\omega (t)} (x)+
r(t,x))$ as in Ansatz (1.3). The reminder $^tR
 =(  r,
  \overline{r} )$ will be split as
$$R=\left ( z_+(t) + {z}_-(t)\sigma _1\right ) \xi (\omega (t), x)+
f_d(t,x)+f_c(t,x)  \tag 4.1$$ with  $f_d +f_c$ the splitting in
(2.8). In analogy to standard constructions of center and stable
manifolds, we consider functional spaces where we will interpret $X$
as the set of fixed points of certain functionals.

   For $p>5$ the exponent in (1.1) and for $4/q=1-1/p$ we set
$$\aligned & \Cal Z :=L_t^4L_x^\infty \cap L^q _t  W^{1,2p}_x   \cap
  L^\infty _t H^1_x\cap C^0_tH^1_x\cap
   H_x^{1,-2} L_t^2
  ([0,\infty )\times \Bbb R, \Bbb C^2);\\&  \widehat{\Cal X} :=\{ ( z_+(t),z_-(t),
  \gamma (t), f_{ c}(t,x)): \,  z_\pm  (t)\in (L^1 \cap
L^\infty \cap C^0)([0,\infty ),\Bbb R) ;\\&  f_{ c}(t,\cdot) \in
L^2_c ( {\omega _0}) \cap \Cal Z\text{ with }
\overline{f_c}=\sigma_1f_c;\,
 \dot \gamma   \in (L^1\cap L^\infty ) ([0,\infty ) ,\Bbb R)
\}    \endaligned $$ with,  for $\widehat{\Cal R}= ( z_+ ,z_- ,
\gamma  , f_{ c} )$,
$$\| \widehat{\Cal R} \| _{\widehat{\Cal X}  }= \| (z _+,z_-)\| _{(L^1 \cap L^\infty )[0,\infty ) }+  \|
  \gamma   \| _{  (W^{1,\infty }\cap \dot  W^{1,1 })
[0,\infty )}   + \| f_{ c} \| _{\Cal Z  }.$$ Let ${\Cal X}
 :=((W^{1,\infty }\cap \dot  W^{1,1 }) [0,\infty ))\times \widehat{\Cal
X}  $ with elements $\Cal R=(\omega , \widehat{\Cal R}).$ Fix
$\alpha \in (0,1)$ and $\omega _0\in (\alpha , 1/\alpha ).$ For
$\epsilon \in (0,\epsilon _0]$ let $$\aligned & B _{\Cal X }
 (\omega _0,\gamma _0,\epsilon )=\{ \Cal R \in \Cal X   :\| \Cal R-(\omega
_0,0,0,\gamma _0,0) \|  _{\Cal X}\le \epsilon \} \\&  {B}
_{\widehat{\Cal X }}
 ( \gamma _0,\epsilon  )=\{ \widehat{\Cal R} \in \widehat{\Cal X}    :\|
\widehat{\Cal R}-( 0,0,\gamma _0,0) \| _{\widehat{\Cal X}}\le
\epsilon \}.\endaligned
$$
For $\epsilon _0<\delta (\alpha)$,  with $ \delta (\alpha)$ chosen
to be the same of  Lemma 2.4, in $B _{\Cal X }
 (\omega _0,\gamma _0,\epsilon )$ by definition we have $\| \omega
(t)-\omega _0\| _\infty < \delta (\alpha)$. By Lemma 2.4 we define
$f_d(t,x)\in L^2_c ( {\omega _0})$ with $|f_d|\ll |f_c|$,  so that
$f(t,x)=f_d(t,x)+f_c(t,x)\in L^2_c ( {\omega ( t)})$ and $ \| f\|
_{\Cal Z   } \approx \|  f_c\| _{\Cal Z  } .$ Then given $\Cal R \in
B _{\Cal X }
 (\omega _0,\gamma _0,\epsilon )$
we define
  $ R(t,x,\Cal R)$  by formula (4.1).
By construction, $R(t )\in N^\perp _g(H _{\omega ( t)}^\ast )$  and
$ \| R   \| _{\Cal Z  } \le C   ( \| (z_+,z_-) \| _{L^1_t\cap
L_t^\infty } + \| f _{ c}  \| _{\Cal Z } ) $ for $C =C (\alpha )$.
We fix $\omega (0)
>0$ (resp.$\gamma
 (0) \in \Bbb R$) close to $\omega _0$ (resp.$\gamma
 _0 $) and
for $R=R(\Cal R)$ we write

$$\align &\omega  (t)=\omega  (0)+\widetilde{\omega}(\Cal R)  , \, \widetilde{\omega} (\Cal R)  := \int
^{t}_0\dot {\widetilde{\omega}} (\Cal R)(s)  ds,
      \tag 4.2 \\& \gamma  (t)=\gamma  (0) +
 \widetilde{\gamma}(\Cal R)\, , \quad \widetilde{\gamma} (\Cal R) := \int _{0}^{t}
 \dot {\widetilde{\gamma}} (\Cal R)(s)  ds ,\tag 4.3
\endalign$$ where $\dot {\widetilde{\omega}} (\Cal R)$
and $\dot {\widetilde{\gamma}} (\Cal R)$ are as in (2.5).
Schematically  we have for $\Cal R \in B _{\Cal X }
 (\omega _0,\gamma _0,\epsilon )$

 $$ \dot {\widetilde{\omega}} (\Cal
R):=\big \langle O(R^2(t )),
 { \Phi _{\omega (t)}  }\big \rangle \text{ and }
 \dot {\widetilde{\gamma}} (\Cal R):=\big  \langle O(R^2(t )),
 { \partial _\omega \Phi _{\omega (t)}  }\big \rangle .\tag 4.4
  $$
  Let $\ell (t,\Cal R)=\omega (t)-\omega _0+\dot {\widetilde{\gamma}} (\Cal
  R)$.
For a small    $h_0  \in H^1 (\Bbb R, \Bbb C^2)\cap L^2_c( {\omega
_0}) $, $h_0(-x)=h_0(x)$,
  with $ \overline{h_0}=\sigma _1 h_0   $,   we write

  $$\aligned & P_\pm (\omega _0)f_{ c}( t,x )=
   e^{-it H_{\omega _0}}e^{\mp i\int _0^t \ell  (\tau ,\Cal R) d\tau
}P_\pm (\omega _0)h_0 (x )- \widetilde{f}_c(\Cal R ),\\&
\widetilde{f}_c(\Cal R ):=  \int _t ^{ \infty }e^{-i(t-s) H_{\omega
_0}} e^{\mp i\int _s^t \ell  (\tau ,\Cal R) d\tau }P_{ \pm}(\omega
_0) F(\Cal R) (s) ds,\\& F(\Cal R):=P_c(\omega _0) \big  \{
P_c(\omega )\sigma _3\dot {\widetilde{\gamma }}(\Cal R) (z _{+}+
z_{-} \sigma _1  ) \xi +  N(R) + i\dot {\widetilde{\omega }}(\Cal R)
\partial _\omega P_c(\omega )f  \\&  +\left ( \omega _{0}
V(\sqrt{\omega _0} x)- \omega V(\sqrt{\omega} x)\right ) f+
\ell(t,\Cal R)  \left ( P_c(\omega )-P_c(\omega _0) \right ) \sigma
_3
 f   \\& + \ell(t,\Cal R)(P_+(\omega
_0)-P_-(\omega _0)-P_c(\omega _0)\sigma _3) f_c\big  \} .
\endaligned \tag 4.5$$
We write
$$ \aligned & z_+ (t)= \widetilde{z}_+(\Cal R),
 \text{ where  }   \widetilde{z}_+ (\Cal R):=\\&
  d_1\int _{t} ^{\infty }ds\, e^{   \int
_{s}^{t}\mu (\omega (s')) ds'} \big \{     i \dot {\widetilde{\omega
}}(\Cal R) \langle f(s ), \sigma _3 \partial _\omega \xi (\omega (s)
) \rangle +  \langle \sigma _3\dot {\widetilde{\gamma }} (\Cal R)
R(s)
  \\& +  N(
R(s))- i \dot {\widetilde{\omega }}(\Cal R) \left [ z_+(s) + z _-(s)
\sigma _1\right ]
\partial _\omega \xi (\omega (s) ) ,
\sigma _3  \xi (\omega (s) ) \rangle \big \} .
\endaligned   \tag 4.6$$
For  a $z_-(0)$ small we  write $$ \aligned & z_-(t)= e^{- \int
_{0}^{t}\mu (\omega (s )) ds }z_{-} (0)+ \widetilde{z}_-(\Cal
R)\text{ where  }   \widetilde{z}_- (\Cal R):=\\&
  d_1\int _{0} ^{t }ds\, e^{-  \int
_{s}^{t}\mu (\omega (s')) ds'} \big \{     i \dot {\widetilde{\omega
}}(\Cal R) \langle f(s ), \sigma _1\sigma _3
\partial _\omega \xi (\omega (s) ) \rangle +  \langle \sigma _3\dot
{\widetilde{\gamma }} (\Cal R)   R(s)
  \\& +  N(
R(s))- i \dot {\widetilde{\omega }}(\Cal R) \left [ z_+(s) + z _-(s)
\sigma _1\right ]
\partial _\omega \xi (\omega (s) ) , \sigma _1\sigma _3 \xi (\omega (s) ) \rangle \big \} .
  \endaligned   \tag 4.7$$
We interpret (4.2-7) as an equation in $  B _{\Cal X}(\omega
_0,\gamma _0,\epsilon )   \subset \Cal X  $.

 \proclaim{Proposition 4.1} Fix $\alpha \in (0,1) $ and $\omega _0\in
 (\alpha , 1/\alpha )$, $\gamma _0\in \Bbb R$.
Then $\exists$ $\epsilon _0 >0$,   $c(\alpha )$  and a $C>0$ such
that $\forall$ $(\omega (0), \gamma (0),z_-(0))$, with $z_-(0)\in
\Bbb R$, $\gamma (0)\in \Bbb R$, with

$$ |\omega (0) -\omega _0|+ |\gamma (0) -\gamma  _0|
+|z_-(0)|<\epsilon /5\le \epsilon _0/5$$
  and $\forall$
$h_0\in  H^1  _r(\Bbb R,\Bbb C^2)\cap L^2_c( {\omega _0 })$
satisfying $ \overline{h_0}=\sigma _1 h_0 $ with $\| h_0 \| _{ H^1
(\Bbb R)}<c(\alpha )\epsilon $, and if we define $R(t,x)$ by (4.1)
with $f_d(t,x)$ defined by Lemma 2.4, then, for $4/q=1-1/p$ with
$p>5$ the exponent in (1.1), there exists a   solution
$$(\omega (t),z_+(t), z_-(t),\gamma (t), f_c(t))\in C^0([0,\infty ),
\Bbb R ^4) \times \Cal Z   \tag 4.8 $$
 of (4.2-7)  such that $\forall$  $t\ge 0$  we have $f_c(t)=\sigma
_1\overline{f_c(t)}$
 and
$$ \align & |  \omega
(t)-\omega _0 |\le  \epsilon  , \,   | (\omega (t), z_-(t), \gamma
(t)) -(\omega (0) , z_-(0), \gamma (0) ) |<C\epsilon ^2 ,\tag 1
\\&
\| f _c\| _{ \Cal Z  }\le \epsilon ,\tag 2\\&  \|  z_-(t)  \| _{(L^1
\cap L^\infty ) [0,\infty )}\le  \epsilon , \, \,  \| z_+(t) \|
_{(L^1 \cap L^\infty ) [0,\infty )}<C \epsilon ^2 \tag 3\\&
 \lim _{t\to
\infty}\left ( z_+(t),z_-(t)   ) =(0,0
  \right ). \tag 4
\endalign$$
There exist $\gamma _\infty \in \Bbb R$, $\omega _\infty >0$ such
that $$\lim _{t\to \infty}(\omega (t),\gamma (t)) = (\omega _\infty
, \gamma _\infty )\tag 5$$ and for $\ell (t )=\omega (t)-\omega
_0+\dot { {\gamma}} (t)$
 $$ \align  \lim _{t\to \infty}\| f(t)- e^{-it
H_{\omega _0}}e^{  i\int _0^t \ell  (\tau ) d\tau (P_- (\omega
_0)-P_+(\omega _0))} h_0 \| _{ H^1(\Bbb R,\Bbb C^2)}=0.\tag 6
\endalign$$
   We have $\sigma _1R(t,x) =\overline{ R(t,x)}  $,
  $R(t,x)$ solves (2.3),  the first entry $r(t,x)$ of $
^tR=(r,\overline{r})$  solves (2.1) and $u(t,x)$, defined in (1.3),
solves (1.1).
\endproclaim
There is an isomorphism between the space of the $u(t,x)$ and the
space of the $U(t,x)={^t(u(t,x), \overline{u}(t,x))}  $, so we think
$X$ in the latter space. The spirit of Proposition 4.1 is that we
try to parametrize the set $X$ by means of $(\omega (0), \gamma (0),
z_-(0),h_0)$. In fact we cannot exclude that for each choice of the
parameter there are more than one solutions of the form (4.8). So we
define the  $X$ in Theorem 1.1 as  the union of the trajectories
associated to all possible solutions of (4.2-7).

\head \S 5 Proof of Proposition 4.1\endhead

  Set $\Cal R
=(\omega , \widehat{\Cal R} ) $ with  $ \widehat{\Cal R} = (z_+
 ,z_-  , \gamma  , f_c  ) $.  Set for $\ell (t,\Cal R)=\omega (0)
 +\widetilde{\omega }(\Cal R)-\omega _0+\dot {\widetilde{\gamma}} (\Cal
  R)$
$$\aligned &
L(\Cal R):= (\widetilde{z}_+(\Cal R),\widetilde{z}_-(\Cal
R),\widetilde{\omega }(\Cal R), \widetilde{\gamma} (\Cal R) );
\\&
\Cal G(\omega )(\widehat{\Cal R}):= ( 0, e^{-\int _0^t\mu (\omega
(s)) ds} z_-(0)    , \gamma (0) , \\& e^{-itH_{\omega _0} }(P_{
+}(\omega _0)e^{- i\int _0^t \ell  (\tau ,\Cal R) d\tau }+P_{
-}(\omega _0)e^{ i\int _0^t \ell  (\tau ,\Cal R) d\tau })h_0)
  +\widetilde{\Cal G} (\Cal R);
\\& \widetilde{\Cal G} (\omega )
(\widehat{\Cal R}):=( \widetilde{z}_+(\omega , \widehat{\Cal R}
),\widetilde{z}_-(\omega , \widehat{\Cal R} ), \widetilde{\gamma}
(\omega , \widehat{\Cal R} ), \widetilde{f}_c(\omega , \widehat{\Cal
R} ));
\\&
\Cal F(\Cal R ):= ( \omega (0), e^{-\int _0^t\mu (\omega (s)) ds}
z_-(0) ,  \gamma (0) , \\& e^{-itH_{\omega _0} }(P_{  +}(\omega
_0)e^{- i\int _0^t \ell  (\tau ,\Cal R) d\tau }+P_{ -}(\omega _0)e^{
i\int _0^t \ell  (\tau ,\Cal R) d\tau })h_0)
  +\widetilde{\Cal F} (\Cal R); \\& \widetilde{\Cal F} (\Cal
R):=(\widetilde{\omega }(\Cal R), \widetilde{\Cal G} (\omega
)(\widehat{\Cal R})).
\endaligned \tag 5.1$$
To prove Proposition 4.1 we look for fixed points of $\Cal F(\Cal R
)$. We are not able to show that $\Cal F(\Cal R )$ is Lipschiz
because of the $\omega (t)$  in the $\ell (t, \Cal R)=\omega (t)
-\omega _0+\dot {\widetilde{\gamma}} (\Cal R)$ and the exponent
$\int _s^t\ell (\tau , \Cal R)d\tau $ in the definition  (4.5) of
$\widetilde{f}_c(\Cal R)$.  We split $\Cal R=(\omega , \widehat{\Cal
R})$  and we   solve the system by substitution, by first solving
for $\widehat{\Cal R}$ with $\omega $ arbitrary but with $\| \omega
- \omega _0\| _\infty   $ small. Since $\Cal F( \omega ,
\widehat{\Cal R})$ is Lipschiz and a contraction in $\widehat{\Cal
R}$, with constant independent of $\omega$, for each $\omega $ we
get a unique corresponding $\widehat{\Cal R}=\widehat{\Cal R}(\omega
)$ by the contraction principle. $ \widehat{\Cal R}(\omega )$ is
continuous in $\omega $. Substituting in the equation for $\omega$,
we obtain a fixed point problem in $\omega$ which we solve by the
Schauder fixed point theorem.

 By Lemmas 3.1-2 we
have: \proclaim{Lemma 5.1} For $\alpha \in (0,1) $ there exists
$C(\alpha )>0$ such that $\forall$ $\omega _0\in (\alpha , 1/\alpha
)$   we have   $  \| e^{-itH_{\omega _0} }P_c(\omega _0)h \| _{\Cal
Z  }< C(\alpha ) \|  h\| _{H^1};$ $e^{-itH_{\omega _0} }$ is
strongly continuous in $H^1(\Bbb R , \Bbb C^2).$
\endproclaim
 Next, we
have:

\proclaim{Lemma 5.2} There exists a fixed $C>0$ such that for all
$0<\epsilon <\epsilon _0$, $  L(\Cal R) $ is $C^1$ in $B _{\Cal
X}(\omega _0,\gamma _0,\epsilon ) $  such that for $L(\Cal R) =
(\widetilde{z}_+(\Cal R),\widetilde{z}_-(\Cal R),\widetilde{\omega
}(\Cal R), \widetilde{\gamma} (\Cal R) )$ and for any $t_0\ge 0$
$$\aligned & \| L(\Cal R)\| _{((L^1\cap L^\infty )^2\times
(W^{1,\infty }\cap W^{1,1})^2)[t_0,\infty ) }\le \\&  C \epsilon (
\epsilon e^{-\frac{\alpha \mu (1)}{2} t_0}+ \| (z_+ , z_- )\|
_{L^1[t_0,\infty )} + \| f _c \| _{L^2((t_0,\infty ) ,
L^{2,-2}_x)}).\endaligned \tag 1$$ Furthermore we have

 $$\aligned &  \| D L(\Cal R)\delta \Cal
R\| _{((L^1\cap L^\infty )^2\times (W^{1,\infty }\cap
W^{1,1})^2)[0,\infty ) }\le C \epsilon \| \delta \Cal R\| _{\Cal X
 } .\endaligned \tag 2$$
\endproclaim

{\it Proof.} Set $\widetilde{z }_+(t) =\widetilde{z }_+(\Cal R)(t)
=d_1\int _t^{+\infty } ds\, e^{ \int _{s}^{t}\mu (\omega (s')) ds'}
Z _+(\Cal R) (s)$,
$$\aligned &
   Z _+(\Cal R ) : =\langle  \sigma _3\dot {\widetilde{\gamma }
}    R +  N (R  )  -i \dot {\widetilde{\omega }}  \left [
 z_+
+   z _-  \sigma _1\right ]  \partial _\omega  \xi (\omega   ) ,
\sigma _3 \xi (\omega  )   \rangle    +i \dot {\widetilde{\omega }}
\langle f , \sigma _3 \partial _\omega \xi (\omega   ) \rangle
.\endaligned
$$
In $B _{\Cal X}(\omega _0,\gamma _0,\epsilon )$  we have $ \mu
(\omega (t))> \alpha \mu (1)>0$. So for $t\ge t_0$
$$\aligned & |\widetilde{z }_+(t)|+ \| \widetilde{z }_+ \| _{L^1[t_0,\infty
)} \le \int _t^{+\infty } ds e^{-\alpha \mu (1)|t-s|}|Z_+(s)|ds+
\frac 1{\alpha \mu (1)} \| Z_+ \| _{L^1[t_0,\infty )} .\endaligned$$
The above is $\le C_\alpha \| Z_+ \| _{L^1[t_0,\infty )} $. We have
$$\| Z_+ \| _{L^1[t_0,\infty )} \le C \epsilon (  \| (z_+ , z_- )\| _{L^1[t_0,\infty )}
+ \| f _c \|  _{L^2((t_0,\infty ) , L^{2,-2}_x)}).
$$
So for $t\ge t_0$ we get $$|\widetilde{z }_+(t)|+ \| \widetilde{z
}_+ \| _{L^1[t_0,\infty )} \le C \epsilon (  \| (z_+ , z_- )\|
_{L^1[t_0,\infty )} + \| f _c \|  _{L^2((t_0,\infty ) ,
L^{2,-2}_x)}).\tag 3 $$ We have $\widetilde{z }_-(t) =\widetilde{z
}_-(\Cal R)(t) =d_1\int _0^{t } ds\, e^{ -\int _{s}^{t}\mu (\omega
(s')) ds'} Z _-(\Cal R) (s)$ with
$$\aligned &
   Z _-(\Cal R )   =\langle  \sigma _3\dot {\widetilde{\gamma }
}    R +  N (R  )  -i \dot {\widetilde{\omega }}  \left [
 z_+
+   z _-  \sigma _1\right ]  \partial _\omega  \xi (\omega   ) ,
\sigma _1\sigma _3 \xi (\omega  )   \rangle    +i \dot
{\widetilde{\omega }} \langle f , \sigma _1\sigma _3 \partial
_\omega \xi (\omega   ) \rangle .\endaligned
$$
 Then

$$\aligned & |\widetilde{z }_-(t)|
\le \int _0^{t } ds e^{-\alpha \mu (1)|t-s|}|Z_-(s)|ds ,\\&  \|
\widetilde{z }_- \| _{L^1[t_0,\infty )}  \le \left \| \int _0^{t }
ds e^{-\alpha \mu (1)|t-s|}|Z_-(s)|ds \right \| _{L^1[t_0,\infty )}
.\endaligned$$ From the first we read for $t\ge t_0$

$$\aligned & |\widetilde{z }_-(t)|
\le  C  e^{-\alpha \mu (1) t /2} \| Z_- \| _{L^1[0,t_0/2 )} + C \|
Z_- \| _{L^1[ t_0/2,\infty )} .\endaligned$$ This yields for $t\ge
t_0$
$$|\widetilde{z }_-(t)|  \le C  \epsilon ^2  e^{-\alpha \mu (1) t /2}+C \epsilon (  \| (z_+ , z_- )\|
_{L^1[t_0,\infty )} + \| f _c \|  _{L^2((t_0,\infty ) ,
L^{2,-2}_x)}).\tag 4 $$ In a similar fashion we obtain
$$  \| \widetilde{z
}_- \| _{L^1[t_0,\infty )}\le C  \epsilon ^2  e^{-\alpha \mu (1) t_0
/2}+C \epsilon (  \| (z_+ , z_- )\| _{L^1[t_0,\infty )} + \| f _c \|
_{L^2((t_0,\infty ) , L^{2,-2}_x )}) . \tag 5$$ Notice that (3-5)
imply
$$\| (z_+ , z_- )\| _{L^1[t_0,\infty )}\le C  \epsilon ^2  e^{-\alpha \mu (1) t_0
/2}+ C  \epsilon \| f _c \| _{L^2((t_0,\infty ) , L^{2,-2}_x)}.\tag
6$$ By (4.4) we have
$$ \aligned & \| \dot {\widetilde{\omega}} (\Cal
R) \| _{L^1[t_0,\infty )} = \| \langle O(R^2(t )),
 { \Phi _{\omega (t)}  }\big \rangle  \| _{L^1[t_0,\infty )} \le C
 \| R\| ^2_{L^2([t_0,\infty ), L^{2,-2}_x)}.
 \endaligned$$
Then
$$   \| \dot {\widetilde{\omega}} (\Cal
R) \| _{L^1[t_0,\infty )}  \le  C \epsilon (  \epsilon e^{-\alpha
\mu (1) t_0 /2} + \| f _c \|  _{L^2((t_0,\infty ) ,
L^{2,-2}_x)}).\tag 7 $$ Similarly

$$   \| \dot {\widetilde{\gamma}} (\Cal
R) \| _{L^1[t_0,\infty )}  \le  C \epsilon (  \epsilon e^{-\alpha
\mu (1) t_0 /2}+ \| f _c \|  _{L^2((t_0,\infty ) ,L^{2,-2}_x)}).\tag
8 $$ Then (3-5) and (7-8) yield  (1).

We have $ \widetilde{z }_+(\Cal R+\delta \Cal R)=\widetilde{z
}_+(\Cal R)+\delta \widetilde{z }_+(\Cal R+\delta \Cal
R)=\widetilde{z }_+(\Cal R)+ D \widetilde{z }_+(\Cal R) \delta \Cal
R+O(\delta \Cal R ^2)   $
$$\aligned & \text{with }D
\widetilde{z }_+(\Cal R) \delta \Cal R= \int _t^{+\infty } ds\, e^{
\int _{s}^{t}\mu (\omega (s')) ds'}\times \\& \left [  \left ( \int
_{s}^{t}  \mu (1) \delta \omega (s') ds'\right )  Z_+ (\Cal R)
(s)+DZ_+ (\Cal R)\delta \Cal R \right ] .\endaligned $$ So $ \| DZ_+
(\Cal R)\delta \Cal R\| _{L^1[0,\infty )}\le \widetilde{C}_\alpha
\epsilon \| \delta \Cal R \| _\Cal X $ and   $ |(D \widetilde{z
}_+(\Cal R)\delta \Cal R )(t)|+ \| D \widetilde{z }(\Cal R) \delta
\Cal R \| _{L^1[0,\infty )}\le \widehat{C}_\alpha \epsilon \| \delta
\Cal R \| _{\Cal X}$. For $ \| \delta \Cal R \| _{\Cal X}\lesssim
\epsilon$, $ |O(\delta \Cal R ^2)(t)|+ \|O(\delta \Cal R ^2) \|
_{L^1[0,\infty )}\lesssim \epsilon \| \delta \Cal R \| _{\Cal X}$.
Similar estimates hold for the $\widetilde{z}_-(\Cal R)$,
$\widetilde{\omega} (\Cal R)$ and $\widetilde{\gamma}(\Cal R)$. This
yields (2).

\bigskip
 Consider the ball $
B_{L^\infty}(\omega _0,\epsilon  )$ defined by
 $\| \omega (t)-\omega _0\| _{ L^{ \infty}[0, \infty )  }<\epsilon  $.
  \proclaim{Lemma 5.3} {\item {(1)}} There is a fixed $C>0$ such
that we have  $ \| \widetilde{f}_c(\Cal R)\| _{\Cal Z } \le
C\epsilon ^2$ for any  $\Cal R\in B _{\Cal X}(\omega _0,\gamma
_0,\epsilon )$. {\item {(2)}} There is a fixed $C>0$ such that given
any $\omega  \in \overline{B_{L^\infty}(\omega _0,\epsilon  )}$ the
map $\widehat{\Cal R} \in B _{\widehat{\Cal X}}( \gamma _0,\epsilon
) \to \widetilde{f}_c(\omega , \widehat{\Cal R})\in \Cal Z $ is
differentiable with   $ \| D\widetilde{f}_c(\omega , \widehat{\Cal
R})\delta \widehat{\Cal R}\| _{\Cal Z } \le C\epsilon
\|\widehat{\Cal R}\| _{\widehat{\Cal X} }$. {\item {(3)}} Let $\Cal
R_j=(\omega , \widehat{\Cal R}_j)  $ with $ \omega \in
B_{L^\infty}(\omega _0,\epsilon  )$ and $  \widehat{\Cal R}_j \in B
_{\widehat{\Cal X}} (\gamma _0,\epsilon ) $ for $j=1,2.$ Then
$$ \| e^{-itH_{\omega _0} } P_{
\pm }(\omega _0)\left ( e^{\mp i\int _0^t \ell  (\tau , {\Cal R}_1)
d\tau }  - e^{\mp i\int _0^t \ell  (\tau ,  {\Cal R}_2) d\tau
}\right ) h_0\| _{\Cal Z}\le C \epsilon\| \widehat{\Cal R}_1-
\widehat{\Cal R}_2\| _{\widehat{\Cal X}} \| h_0\| _{H^1_x}.$$

\endproclaim
{\it Proof.}   (3) follows by

$$ \aligned & \| e^{-itH_{\omega _0} } P_{
\pm }(\omega _0)e^{\mp  i\int _0^t (\omega  (\tau ) -\omega _0)
d\tau } \left ( e^{\mp i\int _0^t   \dot {\widetilde{\gamma} }  (
{\Cal R}_1)(\tau )   d\tau  } -e^{\mp i\int _0^t  \dot
{\widetilde{\gamma} } ( {\Cal R}_2 )(\tau )   d\tau  }\right ) h_0\|
_{\Cal Z}\\& \le C_1 \| \dot {\widetilde{\gamma} }  (  {\Cal
R}_1)-\dot {\widetilde{\gamma} }  (  {\Cal R}_2)   \| _{L^1_t} \|
h_0\| _{H^1_x}\le C _2\epsilon\| \widehat{\Cal R}_1- \widehat{\Cal
R}_2\| _{\widehat{\Cal X}} \| h_0\| _{H^1_x}.
\endaligned $$
The first two claims of Lemma 5.3  are a consequence of Lemmas 5.4
and 5.5 below. We have a decomposition $N(R)=O_{loc}(R^2)+N_2(f_c)$
with $N_2(f_c)=O(f^p_c)$. We set $F(\Cal R)=F_1(\Cal R)+F_2(\Cal R)$
with $F_2(\Cal R)=N_2(f_c)=O(f^p_c)$.

\proclaim{Lemma 5.4} Let $\omega (t)$ be  a   function with values
in $(\alpha , 1/\alpha )$.   Then for a fixed $C=C(\alpha)$ we have
$   \| \widetilde{f} _c(\Cal R)\| _{\Cal Z }
 \le C    \left ( \| F_1(\Cal R) \|  _{    H^{1, 2}_x L_t^2
}+   \| F_2(\Cal R)\|  _{     L^1 _tH^1 _x    } \right ) .$
\endproclaim
{\it Proof.} By Lemmas 3.1, 3.4 and 3.5 for $t_0 \ge 0$
$$ \aligned & \| \widetilde{f}_c (\Cal R)\|
 _{L_t^4((t_0,\infty ),L_x^\infty )\cap L^q _t((t_0,\infty ),
 W^{1,2p}_x )  \cap
  L^\infty _t((t_0,\infty ), H^1_x )   }
 \le \\& \le
 C    \left ( \| F_1(\Cal R)\|  _{  L_t^2(t_0,\infty )   H^{1, 2}_x
}+   \| F_2(\Cal R)\|  _{    L^1 _t((t_0,\infty ),H^1 _x)    }
\right ) .\endaligned
$$ Let $f_j(\Cal R)$ be defined by (4.5) with $F(\Cal R)$ replaced
by $F_j(\Cal R)$. By Lemmas  3.3 and 3.5

 $$ \aligned & \| \widetilde{f} _1\| _{ L_t^2(t_0,\infty ) H_x^{1,-2}
 }\le   C\| F _1\| _{  L_t^2(t_0,\infty )H_x^{1, 2} }\endaligned $$
  By Lemma  3.2,
  for a fixed $C$ and for $t_0\ge 0$,
  $\|  \widetilde{{f}}  _2\| _{L_t^2(t_0,\infty )H_x^{1, -2}}\le $

$$\aligned &   \le \left \|  \int _ t^\infty  ds \|
e^{-i(t-s) H_{\omega _0}} e^{\pm i\int _s^t \ell (\tau , \Cal R)
d\tau} P_\pm (\omega _0 )F_2(\Cal R)(s)\| _{H_x^{1, -2} } \right \|
_{ L_t^2 (t_0,\infty ) }
\\&
\le  \int _{t_0}^\infty ds  \| e^{-i(t-s) H_{\omega _0}} e^{\pm i
\left (\int _0^t  \ell (\tau , \Cal R) d\tau -\int _0^s\ell (\tau ,
\Cal R) d\tau \right) } P_\pm (\omega _0 )F_2(\Cal R)(s)\| _{H_x^{1,
-2}  L_t^2  }
 \\& \lesssim   \int _{t_0}^\infty ds   \| e^{\mp
i\int _{0}^s \ell (\tau , \Cal R) d\tau}F_2(\Cal R)(s)\| _{H^1_x} ds
=\| F_2(\Cal R)\| _{ L^1 _t((t_0,\infty ),H^1 _x) }.
\endaligned $$

The final step to prove Lemma 5.3 is: \proclaim{Lemma 5.5} The maps
$F_j(\Cal R) $ are for $\Cal R\in B _{\Cal X }
 (\omega _0,\gamma _0,\epsilon )$ continuous and
 differentiable, with target
$L_t^2 H_x^{1,  2}$ for $F_1(\Cal R) $ and $L^1 _tH^1 _x$ for
$F_1(\Cal R) $. There exists $C>0$ such that for $\Cal R\in B _{\Cal
X }
 (\omega _0,\gamma _0,\epsilon )$ we have for $t_0\ge 0$, $p>5$
  the exponent in (1.1) and for $4/q=1-1/p$

 $$ \align  & \| F_1(\Cal R) \|  _{    L_t^2(t_0,\infty )H_x^{1,  2}} \le
 C \epsilon (
\epsilon e^{-\frac{\alpha \mu (1)}{2} t_0}+  \| f _c \|
_{L_t^2(t_0,\infty )L_x^{2, -2})    }) \tag 1\\&
   \|
F_2(\Cal R)\|  _{      L^1 _t((t_0,\infty ),H^1 _x )   }    \le
 C \epsilon   \| f _c \| _{L^q((t_0,\infty ) ,  W^{1,2p}_x)    }).
\tag 2 \endalign
$$ We have
$$\align \| DF_1(\Cal R)\delta R\|
_{      L_t^2[0,\infty )H_x^{1,  2} }+ \| DF_2(\Cal R)\delta R\| _{
L^1 _t([0,\infty ),H^1 _x) } \le C\epsilon  \| \delta R\| _{\Cal X }
.\tag 3\endalign
$$
\endproclaim
{\it Proof.} By Lemma 5.12 \cite{C1}, repeated in Appendix B in
\cite{C3}, for $C_{M,N}(\omega )$ upper semicontinuous in $\omega $,
$\forall $ $M$ and $N$ we have
$$ \| \langle x \rangle ^{N} (P_+(\omega  )-P_-(\omega  )-P_c(\omega
)\sigma _3) f\| _{L^2_x}\le C_{M,N}(\omega _1) \|  \langle x \rangle
^{-M} f\| _{L^2_x} . $$ Schematically we have $ F_1(\Cal R)
=O(\epsilon ) \psi f+O_{loc}(\Cal R^2) $ for an exponentially
decreasing
  $\psi (x)$.
Then by Lemma 5.2 we get (1). We have $F_2(\Cal R)= O(f^p_c)$ and
this yields

 $$ \|      F_2  \| _{  L_t^1H^1_x }\lesssim
 \|      f _{c}^p \| _{  L_t^1H^1_x }
\lesssim \left \| \| f_ c\| _{W^{1,2p}_x}    \| f_c\| ^{p-1}_{L^{
2p}_x} \right \| _{L^1_t } \le \| f_c\| _{L^q_tW^{1,2p}_x} \| f_c\|
^{p-1}_{L _t ^{q'(p-1)}  L^{ 2p}_x}.$$ Since $q=\frac{4p}{p-1}
<\frac{4p(p-1)}{3p+1}=q'(p-1)$ by $p>5$, then for some $0<\vartheta
<1 $ we get $ \|      F_2  \| _{ L_t^1H^1_x }\lesssim \| f_c\| ^{1+
\vartheta (p-1)}_{L^q_tW^{1,2p}_x} \| f_c\| ^{(1-\vartheta )
(p-1)}_{L^\infty_tH^{1 }_x}$.  This yields (2). Proceeding similarly
we get (3).

 \proclaim {Lemma 5.6}  Consider $\Cal G (\omega  )$ defined
by (5.1). {\item {(1)}} $\forall$
 $\omega  \in \overline{B_{L^\infty}(\omega _0,\epsilon  )}$
   $\exists$
$\widehat{\Cal R} (\omega   )= (z_+ (\omega ),z_- (\omega ), \gamma
(\omega ), f_{c,\omega} ) \in \widehat{\Cal X}  $, unique,
 such that
$ \widehat{\Cal R} (\omega ,h_0)  \in B _{\widehat{\Cal X} }  (
\gamma _0 ,\epsilon /2) $ satisfies the fixed point problem $
\widehat{\Cal R} (\omega )= \Cal G (\omega  ) \left ( \widehat{\Cal
R} (\omega   ) \right ).$

{\item {(2)}}The map $\omega \in \overline{B_{L^\infty}(\omega
_0,\epsilon
  )}\to \widehat{\Cal R}(\omega
)\in {B} _{\widehat{\Cal X }}
 ( \gamma _0,\epsilon  )$ is continuous.
\endproclaim
{\it Proof.}  For $\epsilon \in (0, \epsilon _0)$ with $\epsilon
_0>0$ small enough, $  {\Cal G} (\omega  )$ maps $ B _{\widehat{\Cal
X} }  ( \gamma _0 ,\epsilon /2) $ into itself. By the estimates on
the derivatives   in Lemmas 5.2 and 5.5, $ \| {\Cal G} (\omega  )
\widehat{\Cal R}_1-{\Cal G} (\omega  ) \widehat{\Cal R}_2\|
_{\widehat{\Cal X} }\le C \epsilon \|
 \widehat{\Cal R}_1-
\widehat{\Cal R}_2\| _{\widehat{\Cal X} }. $ There is a fixed point,
which we denote by $\widehat{\Cal R} (\omega  )$, and which is
unique. This yields (1). Let $C \epsilon <1/2$. We have $ \|
\widehat{\Cal R} (\omega _1 )-\widehat{\Cal R} ( {\omega} _2)\|
_{\widehat{\Cal X} }\le$
$$\aligned & \le \|   {\Cal G} (\omega _1)
\widehat{\Cal R} (\omega _1)- {\Cal G} ( {\omega}_2 ) \widehat{\Cal
R} (\omega _1)\| _{ \widehat{\Cal X} } + \| {\Cal G} ( {\omega} _2)
\widehat{\Cal R} (\omega _1)- {\Cal G} ( {\omega}_2 ) \widehat{\Cal
R} ( {\omega} _2)\| _{\widehat{\Cal X} }\\& \le \| {\Cal G} (\omega
_1) \widehat{\Cal R} (\omega _1)- {\Cal G} ( {\omega} _2)
\widehat{\Cal R} (\omega _1)\| _{\widehat{\Cal X} } +C \epsilon \|
\widehat{\Cal R} (\omega _1)-\widehat{\Cal R} ( {\omega} _2)\|
_{\widehat{\Cal X} }.\endaligned
$$
To complete Lemma 5.6 we need to show that $\omega \in
\overline{B_{L^\infty}(\omega _0,\epsilon
  )}\to \Cal G
(\omega ) \widehat{\Cal R} _0\in  \widehat{\Cal X}$ is continuous
for fixed $\widehat{\Cal R} _0$. In view of Lemma  5.2   it remains
to show the following:

\proclaim {Lemma 5.7} The  map $\Cal R\in B _{\Cal X}(\omega
_0,\gamma _0,\epsilon ) \to \widetilde{f}_c(\Cal R)\in \Cal Z$ is
continuous.
\endproclaim
{\it Proof.} We   write $\Cal R =(\omega , \widehat{\Cal R} ) $ to
distinguish between $\omega $ and $ \widehat{\Cal R}= (z_+
 ,z_-  , \gamma  , f_c  ) $.
By Lemma 5.5, to complete the proof of the continuity of
$\widetilde{f}(\Cal R)$ it is enough to show that for fixed $\Cal
R_0 =(\omega _0, \widehat{\Cal R} _0)$ and if we set   $\Cal R_1
=(\omega _0+\delta \omega , \widehat{\Cal R} _0)$, for any
$\varepsilon >0$ there is $\delta >0$ such that $
|\widetilde{f}(\Cal R _0)-\widetilde{f}(\Cal R_1)| \le \varepsilon $
  if $\| \delta \omega \| _{L^\infty}<\delta  $.
For $g(s)= e^{\mp    i\int _0^s  \delta \omega
  (\tau ) d\tau }P_{ \pm}(\omega _0)F(\Cal R_0)(s)$ we need to
  show that  for any $\varepsilon
>0$ there exists $\delta
>0$ such that $\| \delta \omega \| _{L^\infty}<\delta  $ implies
$$\left \| \int _t^\infty e^{-i(t-s) H_{\omega _0}}  \left ( e^{\pm i\int _s^t  \delta \omega
  (\tau ) d\tau } -1\right ) g(s) ds\right \| _{\Cal Z }<\varepsilon
  .$$ We   fix a large
 number $M>0$.    Then, for $\delta  >0$ with $ M\delta    C \|
F(\Cal R_0)\| _{ \Cal Z}\le \varepsilon  /2 $ and since $ \| g\| _{
H_x^{1, - 2}L^2_t(I)+L^1_t(I,H^1_x)}\le \|  F(\Cal R_0)\| _ {
H_x^{1,   2}L^2_t(I)+L^1_t(I,H^1_x)}  $ for any interval $I$,  we
conclude
$$ \left \|    \int  _{t }^{  t+M
}  e^{-i(t-s) H_{\omega _0}}  \left ( e^{\pm i\int _s^t \delta
\omega
  (\tau ) d\tau } -1\right ) g(s) ds\right \| _{\Cal Z}
  <\varepsilon  /2.$$ We have
$$\aligned & \left \|  \int  _{ t+M}^{\infty }
  e^{-i(t-s) H_{\omega _0}}  \left ( e^{\pm i\int _s^t \delta  \omega
  (\tau ) d\tau } -1\right ) g(s) ds\right \| _{\Cal Z}\\&  \le C   \| F (\Cal R_0)    \|
_{H_x^{1,   2}L^2_t (M,\infty ) + L^1_t((M,\infty ), H^1_x) } \to 0
\text{ for $M\nearrow \infty$.}
  \endaligned  $$
\bigskip
Having   $\widehat{\Cal R} (\omega ) $ for any $\omega \in
\overline{B_{L^\infty}(\omega _0,\epsilon  )}$  we substitute
$\widehat{\Cal R}=\widehat{\Cal R} (\omega ) $ in the system and we
reduce to a fixed point problem in $\omega$.  We  will denote by
$\Cal Z (t_0)$ the space defined like $\Cal Z$ in \S 4 but with the
time interval $(0,\infty)$ replaced by $(t_{0},\infty)$. We get:
\proclaim {Lemma 5.8} There is a $\omega (t)\in
 {B_{L^\infty}(\omega _0,\epsilon /2)}$  such that, for $\Cal
R =(\omega ,\widehat{\Cal R}(\omega
 ))$ and $R(t)=R(t,x,\Cal R)$,
 $$   \omega  (t)= \omega  (0)  + {\widetilde{\omega }}( {\Cal R} )(t)
 \, , \quad    {\widetilde{\omega }}( {\Cal R} )(t)
 =\int _0^t
  \big  \langle O(R^2(s )),
 { \Phi _{\omega (s)}  }\big \rangle    ds .   \tag 1$$
 \endproclaim

{\it Proof.} The map on the right side in (1) sends
$\overline{B_{L^\infty}(\omega _0,\epsilon  )}$ into itself.
  Lemma 5.8 is a consequence of the      Schauder fixed point theorem
  if we are able to show that the image of
  $\overline{B_{L^\infty}(\omega _0,\epsilon  )}$, which we denote
  by $A$, has compact closure  in  $\overline{B_{L^\infty}(\omega _0,\epsilon  )}$.
   First of all, $ A\subset
\overline{B_{L^\infty}(\omega _0,\epsilon /3)}\cap
 (W^{1,\infty}\cap \dot W^{1,1})  $.
  It will be
enough to show that, for any $\varepsilon >0$ there exists
$t_0=t_0(\varepsilon )$ such that for any $\omega \in A$ we have $\|
\dot \omega \| _{L^1_t(t_0,\infty )}<\varepsilon .$ This reduces to
showing that for any $\varepsilon
>0$ there is $t_0>0$ such that
for any $\omega \in  {B_{L^\infty}(\omega _0,\epsilon  )}$, given
  the corresponding $\Cal R=(\omega , \widehat{\Cal R} (\omega ))$,
we have $ \| \widetilde{f}_c( {\Cal R} ) \| _{\Cal Z (t_0)
}<\varepsilon .$  But by the proof of Lemma 5.4 and by (1-2) Lemma
5.5 we get

$$\aligned & \| \widetilde{f}_c (\Cal R)\|
 _{\Cal Z ( t_0)  }
 \le   C    \left ( \| F_1(\Cal R)\|  _{
   H_x^{1,   2} L_t^2(t_0,\infty )
}+   \| F_2(\Cal R)\|  _{    L^1 _t( (t_0,\infty ),H^1 _x )   }
\right ) \le \\& \le  C \epsilon  (\epsilon  e^{-\frac{\alpha t_0
\mu (1)}{2}}+   \| e^{-iH_{\omega _0}t}h_0\|
 _{\Cal Z (t_0)  } +   \| \widetilde{f}_c (\Cal R)\|
 _{\Cal Z (t_0)  })  \endaligned
$$
which implies $\| \widetilde{f}_c (\Cal R)\|
 _{\Cal Z (t_0)  }
 \le   C_1 \epsilon  (\epsilon  e^{-\frac{\alpha t_0
\mu (1)}{2}}+   \| e^{-iH_{\omega _0}t}h_0\|
 _{\Cal Z (t_0)  } ) $
and yields the desired result.

\bigskip By Lemmas 5.6-8 we conclude that we have a solution $\Cal
R=(\omega ,\widehat{\Cal R})\in B _{ {\Cal X} }  ( \omega _0,\gamma
_0 ,\epsilon    )$ which yields the solution (4.8) of Proposition
4.1. Estimates (1-4) as well as the limits (5)  follow from the
definition of $\Cal X$.  Now we prove the remaining part of
Proposition 4.1.   We can define a smooth diffeomorphism from  a
 neighborhood of $(\omega _0,0,0, \gamma _0, 0)\in \Bbb R^4\times
 (H^1_r(\Bbb R ,\Bbb C^2)\cap L^2(\omega _0))$
 with values in  a small neighborhood of  $e^{i\gamma _0}\phi _{\omega
 _0}(x)\in H^1_r(\Bbb R, \Bbb C)$ which associates to every $\Pi = (\omega ^{ (0) } ,z_+^{(0)}, z_-^{(0)},
 \gamma ^{ (0) } , f_c^{(0)}(x)) $
 $$u_\Pi (x)= e^{i\gamma   ^{(0)}} \left ( \phi _{\omega
  ^{(0)}}(x) +r_\Pi ( x) \right )$$ with $^t(r_\Pi ( x),\overline{r}_\Pi (
x))= R_\Pi ( x)$ and, for $f_d[\Pi ](x)$ defined by Lemma 2.4, with
$$R_\Pi ( x)= \left ( z_+^{(0)}  + z_-^{(0)} \sigma _1\right )  \xi (\omega   ^{(0)},x)
+f_d[\Pi ](x)+f_c ^{(0)}(x)  .$$ Then given the solution in (4.8)
and given  $R(t)$ defined by (4.1), the corresponding point in
$u(t)\in H^1_r (\Bbb R, \Bbb C)$ is given by $^t(u,\overline{u})=
e^{i\sigma _3(\int _0^t \omega (s)ds+\gamma (t))}  ( \Phi _{\omega
(t)}+ R(t))$. In particular $u(t)\in C^0([0,\infty), H^1_r (\Bbb R,
\Bbb C))$ and is the solution of (1.1) with $u(0)=u_\Pi $.  By
construction

$$\lim _{t\to  \infty }
\| R(t)  - e^{-it H_{\omega _0}}e^{- i\int _0^t \ell  (\tau ) d\tau
(P_+ (\omega _0)-P_- (\omega _0)) }  h_0\| _{H^1 (\Bbb R, \Bbb C
^2)} =0.$$ For $h_0 =W(\omega _0)\widetilde{h}_0$ with $ W(\omega
_0)=strong-\lim _{t\to  \infty }e^{ itH_{\omega _0}} e^{-it \sigma
_3(-\partial ^2_x+\omega _0)} $, see \cite{C1},  $$ \lim _{t\to
 \infty }  \| f_c(t)  -e^{- i(\int _0^t \omega   (\tau ) d\tau
+\gamma (t)-\gamma (0)-t\omega _0) \sigma _3 }e^{ it \sigma _3(
\partial ^2_x-\omega _0)} \widetilde{h}_0\| _{H^1 (\Bbb R, \Bbb C ^2)}  =0.$$
So for $^t(r_\infty , \overline{r}_\infty )=e^{i\gamma (0) \sigma
_3} \widetilde{h}_0$ and $ ^t(r  , \overline{r} )= R$  we conclude

$$ \lim _{t\to  \infty }  \| e^{ i \int
_0^t \omega   (\tau ) d\tau  +i\gamma (t)}r(t)  -e^{ it
\partial ^2_x } r_\infty\| _{H^1 (\Bbb R, \Bbb C  )}  =0.$$

\head Errata in paper \cite{C1}\endhead

Unfortunately paper \cite{C1} has many mistakes. Fortunately all of
them can be corrected. Among the various mistakes   we list:

{\item{(1)}}Various formulas between sections 5 and 8 are wrong, for
example the formula for the Wronskian from \S 5 on.

{\item{(2)}} In formula (8.2) in \cite{C1} there is a missing term
on the right hand side.

{\item{(3)}} The  really serious mistake is Lemma 5.4 \cite{C1}: not
only the proof is incorrect, but probably the statement is
incorrect.

\bigskip
In \cite{C3} we have revised \cite{C1} simplifying considerably the
argument. In particular the  smoothing   estimates in \S 3
\cite{C1}, which are analogues of   estimates in \cite{M}, have been
replaced by weaker  estimates estimates in \S 3 \cite{C3}. These new
estimates are listed in \S 3 in the present paper and are simple  to
prove. The estimates in \S 3 \cite{C3} are sufficient for the main
result in \cite{C1,C3}. In particular in \cite{C3} most of the
material in   sections from 5 to 8   in \cite{C1} is eliminated. In
particular the statements in \S 3 \cite{C3} are proved immediately
in \S 3 \cite{C3} with elementary arguments based on material
already in the literature. \cite{C3} relies more on \cite{KS}. The
statement that the linear part in \cite{C1} is proven also when the
matrix potential $V(x)=H_\omega - \sigma _3(-\partial _x^2+\omega )$
is not necessarily even, does not stand any more, since \cite{KS}
assumes symmetry of $V(x)$ as an hypothesis. In fact the arguments
from \S 5 to \S 8 in \cite{C1} can be saved in a corrected form, and
this is done in \cite{CV}. However in the present paper we assume
the results in
 \cite{C3}.

\enddocument